\documentclass{amsbook}%
\usepackage{amsmath}%
\usepackage{amsfonts}%
\usepackage{amssymb}%
\usepackage{graphicx}
\theoremstyle{plain}

\newtheorem{definition}{Definition}

\newtheorem{lemma}{Lemma}

\newtheorem{theorem}{Theorem}
\numberwithin{equation}{section}
\begin{document}
\frontmatter
\title[Short Title]{Methods Of Solving Ill-Posed Problems}
\author[S. Srinivasamurthy]{Suresh B Srinivasamurthy\\
B.S., Bangalore University, 1984\\
M.S., Bangalore University, 1986\\
M.Phil., Bangalore University, 1996\\
A REPORT\\
Submitted in partial fulfillment of the\\
requirements for the degree\\
MASTER OF SCIENCE\\
Department of Mathematics\\
College of Arts and Sciences\\
Kansas State University \\
Manhattan, Kansas\\
2004\\
Approved by\\
Major Professor\\
Alexander G.Ramm}
\email[S. Srinivasamurthy]{suresh@math.ksu.edu}
\urladdr{http://www.math.ksu.edu/~suresh}
\subjclass{47A52, 47B25, 65R30}
\keywords {Ill-posed Problems, Variational Regularization, Quasi-solution, Iterative Regularization, Dynamical systems method (DSM), Regularization parameter, Discrepancy principle, operator equations}

\begin{abstract}
Many physical problems can be formulated as operator equations of the form $Au=f$.  If these operator equations are ill-posed, we then resort to finding the approximate solutions numerically.  Ill-posed problems can be found in the fields of mathematical analysis, mathematical physics, geophysics, medicine, tomography, technology and ecology.  The theory of ill-posed problems was developed in the 1960's by several mathematicians, mostly Soviet and American.  In this report we review the methods of solving ill-posed problems and recent developments in this field. We review the variational regularization method, the method of
quasi-solution, iterative regularization method and the dynamical systems method.  We focus mainly on the dynamical systems method as it is found that the dynamical systems method is more efficient than the regularization procedure.
\end{abstract}
\maketitle
\tableofcontents

\chapter*{Acknowledgements}

\markboth{ACKNOWLEDGEMENTS}{ACKNOWLEDGEMENTS}    It gives me immense pleasure to express my deepest sense of
gratitude and indebtedness to my advisor Professor Alexander G.
Ramm for his wise and thoughtful guidance and for the excellent,
motivating lectures he gave, otherwise this report and my
knowledge and appreciation of this subject wouldn't have been
possible.   I offer my sincere regards and thanks to him.\\

I am deeply obliged and greatly indebted to my family and to all
those for their sustained co-operation and inspiration.  There are
no words to express how much I am grateful to them.\\

Finally, I would like to express my profound appreciation and
thanks to the faculty and staff of the Department of Mathematics
and to faculty and to the staff of Kansas State University who
have helped me directly and indirectly in the preparation of the
report and
pursuit of knowledge.\\

Manhattan

July 2004
\qquad\qquad\qquad\qquad\qquad\qquad\qquad\qquad\qquad
Suresh.B.Srinivasamurthy
\chapter*{Dedication}
\emph{\begin{center}DEDICATED WITH LOVE TO\end{center}}
\emph{\begin{center}SNEHA AND VANI\end{center}}
\mainmatter


\chapter{What are ill-posed problems?}
\pagenumbering{arabic}
\section{Introduction}\label{S:intro}
\begin{definition}\label{D1:ill}
A problem of solving an operator equation
\begin{equation}\label{E1:oper}
Au=f
\end{equation}
where $A:X{\longrightarrow Y}$ is an operator from a Banach space
$X$ into a Banach space $Y$ is called well-posed in the sense of
J.Hadamard (1902) iff it satisfies the following conditions.
\begin{enumerate}
\item (\ref{E1:oper}) is solvable for any  $f \in Y$.
(Existence of the solution, i.e.,
$A$ is surjective.)
\item The solution to (\ref{E1:oper}) is unique.
($A$ is injective.)
\item (\ref{E1:oper}) is stable with respect to small perturbations of $f$.
(Continuous dependence of the solution i.e., $A^{-1}$ is
continuous.)
\end{enumerate}
If any of these conditions fails to hold, then problem
(\ref{E1:oper}) is called ill-posed.
\end{definition}

If the solution does not depend continuously on the data, then
small errors, whether round off errors, or  measurement errors, or
perturbations caused by noise, can create large deviations in the
solutions.  Therefore the numerical treatment of ill-posed
problems is a challenge.  We shall briefly discuss below some of
the concepts and auxiliary results used in this report.

Henceforth $D(A),$ $R(A),$ and $N(A):=\{u:Au=0\}$ denote the
domain, range and null-space of $A$ respectively.  Let
$A^{*}:Y^{*}\longrightarrow{X^{*}}$ be the adjoint operator.  For
simplicity we assume below that $X=Y=H$, where $H$ is a Hilbert
space.  If $A$ is \emph{self-adjoint} then $A=A^{*}$.  If $A$ is
injective, then $N(A)=\{0\}$, $A^{-1}$ is well-defined on $R(A)$
and $u=A^{-1}f$ is the unique solution to equation (\ref{E1:oper})
for $f\in{R(A)}$.  Equation (\ref{E1:oper}) is called
\emph{normally solvable} iff $R(A)=\overline{R(A)}$ i.e., iff
$f\perp{N(A^{*})}={\overline{R(A)}}^{\perp}$.  The overbar denotes
closure.  If $N(A)\neq\{0\}$ define the \emph{normal solution
(pseudo-solution)} $u_{0}$ to equation (\ref{E1:oper}) as the
solution orthogonal to $N(A)$.  Then the normal solution is unique
and has the property that its norm is minimal:
$min\|u\|=\|u_{0}\|$, where the minimum is taken over the set of
all solutions to equation (\ref{E1:oper}).  The normal solution to
the equation $Au=f$ can be defined as the least squares solution:
$\|Au-f\|=min.,u\perp{N(A)}$.  This solution exists, is unique and
depends continuously on $f$, if $H$ is finite-dimensional. The
normal solution is also called \emph{minimal-norm solution}.

$A$ is called \emph{closed} if $\{u_{n}\longrightarrow{u},
Au_{n}\longrightarrow{f}\}$ implies
$\{u\in{D(A)\quad{and}\quad{Au=f}}\}$.  By Banach theorem, if $A$
is a linear closed operator defined on all of $X$ then $A$ is
bounded.  A is called \emph{compact} if it maps bounded sets into
pre-compact sets.  The set $\{f\}$ is \emph{bounded} means there
exists $\rho{}>0$ such that $\|f\|\leq\rho$ and a set is
\emph{pre-compact} if any subsequence from the set contains a
convergent subsequence. In a finite-dimensional Banach space a set
is pre-compact iff it is bounded.  If $A$ is an injective linear
compact operator on an infinite dimensional space then $A^{-1}$ is
unbounded.

\emph{Singular Value Decomposition}:  Suppose $A$ is a compact
linear operator on $H$, then $A^{*}$ is compact and
$|A|:=[A^{*}A]^{1/2}$ is self-adjoint, compact and non-negative
definite, $|A|\phi_{j}=\lambda_{j}\phi_{j}$,
$\lambda_{1}\geq\lambda_{2}\geq{...}\geq{0}\rightarrow{0}$ are the
eigenvalues of $|A|$ with \emph{s-values of $A$}:
$s_{j}=s_{j}(A):=\lambda_{j}(|A|)$ and  $\phi_{j}$ are the
normalized eigenvectors of $|A|$. The faster the s-values go to
zero the more  ill-posed problem (\ref{E1:oper}) is. Any bounded
linear operator $A$ admits the polar representation $A:=U|A|$,
where $U$ is an isometry from $R(A^{*})$ onto
$R(A),\|Uf\|=\|f\|,\|U\|=1$.  One has
$|A|=\sum_{j=1}^{\infty}s_{j}(.,\phi_{j})\phi_{j}$; then the
\emph{SVD of A} is
$A=U|A|=\sum_{j=1}^{\infty}s_{j}(.,\phi_{j})\psi_{j}$, where
$\psi_{j}:=U\phi_{j}$.

A closed set $K$ of X is called a \emph{compactum} if any infinite
sequence of its elements contains a convergent subsequence.  A
sequence $u_{n}$ in $U$ \emph{converges weakly} to $u$ in $U$ iff
$\lim_{n\rightarrow{\infty}}(u_{n},\phi)=(u,\phi)$ for all
$\phi{\in{H}}$.  We denote the weak convergence by
$u_{n}\rightharpoonup{u}$.  If $u_{n}\rightharpoonup{u}$ then
$\|u\|\leq\liminf_{n\rightarrow\infty}\|u_{n}\|$.  If
$u_{n}\rightharpoonup{u}$ and $A$ is a bounded linear operator,
then $Au_{n}\rightharpoonup{Au}$.  A bounded set in a Hilbert
space contains a weakly convergent subsequence.  A functional
$F:U\rightarrow{\mathbb{R}}$ is called \emph{convex} if the domain
of $F$ is a linear set and for all $u,v\in{D(F)}$,
$F(\lambda{u}+(1-{\lambda})v)\leq{\lambda{F(u)}}+(1-\lambda{)}F(v)$,
$0\leq\lambda\leq{1}$.

A functional $F(u)$ is called \emph{weakly lower semicontinuous}
from below in a Hilbert space if $u_{n}\rightharpoonup{u}$ implies
$F(u)\leq\liminf_{n\rightarrow{\infty}}F(u_{n})$.  A functional
$F:U\rightarrow{\mathbb{R}}$ is \emph{strictly convex} if
$F(\frac{u_{1}+u_{2}}{2})<\frac{F(u_{1})+F(u_{2})}{2}$ for all
$u_{1},u_{2}\in{D(F)},$ provided $u_{1}\neq{\lambda{u_{2}}}$,
$\lambda=$constant.

Let $F:X\rightarrow{Y}$ be a functional.  Suppose
$F'(u)\eta=\lim_{\epsilon{\rightarrow}+0}\frac{F(u+\epsilon{\eta})-F(u)}{\epsilon}$,
exists for all $\eta$, and $F'(u)$ is a linear bounded operator in
$H$.  Then $F'(u)$ is called \emph{Gateaux} derivative.  It is
called \emph{Fr\'echet} derivative if the limit is attained
uniformly with respect to $\eta$ running through the unit sphere.

\emph{Spectral Theory}: Let $A$ be a self-adjoint operator in a
Hilbert space $H$.  To $A$ there corresponds a family
$E_{\lambda}$ of ortho-projection operators such that
$\phi{(A)}:=\int_{-\infty}^{\infty}\phi{({\lambda})}dE_{\lambda}$;
$\phi{(A)f}:=\int_{-\infty}^{\infty}\phi{({\lambda})}dE_{\lambda}f$;
$D(\phi{(A)})=\{f:\|\phi{(A)}f\|^{2}=\int_{-\infty}^{\infty}|\phi{({\lambda})}|^{2}(dE_{\lambda}f,f)<{\infty}\}$;
$\|\phi{(A)}\|=\sup_{\lambda}|\phi{({\lambda})}|$.

In particular $A=\int_{-\infty}^{\infty}\lambda{dE_{\lambda}}$ and
$D(A)=\{f:\|Af\|^{2}=\int_{-\infty}^{\infty}|{\lambda}|^{2}(dE_{\lambda}f,f)<\infty\}$.
$E_{\lambda}$ is called the \emph{resolution of the identity}
corresponding to the self-adjoint operator $A$.  $\lambda$ is
taken over the spectrum of $A$.  We have
$\rho{({\lambda})}=(E_{\lambda}f,f)$;
$(dE_{\lambda}f,f)=d(E_{\lambda}f,f)=d\rho{({\lambda})}$;
$I=\int_{-\infty}^{\infty}dE_{\lambda}$;
$f=\int_{-\infty}^{\infty}dE_{\lambda}f$;
$\|f\|^{2}=\int_{0}^{\|A\|}d(E_{\lambda}f,f)$; $E_{+\infty}=I$;
$E_{-\infty}=0$; $E_{{\lambda}-0}=E_{\lambda}$.

Let the operator equation $Au=f$ be solvable (possibly
non-uniquely).  Let $y$ be its minimal-norm solution,
${y}\perp{N(A)}$.  Let $B=A^{*}A\geq{0}$ and $q:=A^{*}f$.  Then
$Bu=q$.  Also,

$(B+\alpha{)}^{-1}:=(B+\alpha{I)}^{-1}$ is a positive definite
operator and is given by

$(B+\alpha{)}^{-1}=\int_{0}^{\infty}\frac{dE_{\lambda}}{{\lambda}+\alpha}$.

$\|(B+{\alpha})^{-1}\|=\sup_{0\leq\lambda\leq\|B\|}|\frac{1}{\lambda{+\alpha}}|\leq{\frac{1}{\alpha}},{\alpha}>0.$

$\|[(A+{\alpha})^{-1}A-I]f\|^{2}={\alpha}^{2}\|(A+{\alpha})^{-1}f\|^{2}$
$={\alpha}^{2}\int_{0}^{\|A\|}\frac{d(E_{\lambda}f,f)}{({{\lambda}+\alpha})^{2}}$
$\leq\int_{0}^{\|A\|}{d(E_{\lambda}f,f)}=\|f\|^{2}<\infty$.

If $\alpha\rightarrow{0}$, the integrand
$\frac{{\alpha}^{2}}{(\lambda{+\alpha})^{2}}$ tends to zero.  So,
by the Lebesgue dominated convergence theorem
$\|[(A+{\alpha})^{-1}A-I]f\|^{2}\rightarrow{0}$ as
$\delta\rightarrow{0}$, provided that
$\int_{0}^{0^{+}}d(E_{\lambda}f,f)=0$, i.e., $f\perp{N(A)}$.
\begin{lemma}
If the equation $Au=f$ is solvable then it is equivalent to the
equation $Bu=q$.
\end{lemma}
\emph{Proof}: $({\Rightarrow})Au=f$, so $A^{*}Au=A^{*}f$.

$({\Leftarrow})A^{*}Au=A^{*}Ay$, so $A^{*}A(u-y)=0$, hence
$A(u-y)=0$ and hence $Au=Ay=f$.  Thus we have proved the lemma.
$\Box$

The mapping $A^{+}:f\longrightarrow{u_{0}}$ is called the
\emph{pseudo-inverse} of $A$.  $A^{+}$ is a bounded linear
operator iff it is normally solvable and $R(A)$ is closed.  So
equation (\ref{E1:oper}) is ill-posed iff $A^{+}$ is unbounded.
One can find the details of this in \cite{aR04a}.

An operator $\mathbf{\Phi{(t,u)}}$ is \emph{locally Lipschitz}
with respect to $u\in{H}$ in the sense
$sup\|\mathbf{\Phi(t,u)}-\mathbf{\Phi(t,v)}\|\leq{c}\|u-v\|,
 c=c(R,u_{0},T)>0$ where the supremum is taken for all
$u,v\in{B(u_{0},R)},$ and $t\in{[0,T]}$.
\newpage

\section{Examples of ill-posed problems}\label{Ex1:eg}
\emph{Example 1.  Stable numerical differentiation of noisy
data}

The problem of numerical differentiation is ill-posed in the sense
that small perturbations of the function to be differentiated may
lead to large errors in its  derivative.  Let $f\in{C^{1}}[0,1]$,
with noisy data $\{{\delta}, f_{\delta}\}$, where ${\delta}>0$ is
the noise level, that is we have the estimate
\begin{equation}\label{E2:fnois}
\|f_{\delta} - f\|\leq{\delta}.
\end{equation}
The problem is to estimate stably the derivative $f^{'}$, i.e., to
find such an operation $R_{\delta}$ such that the error estimate\\
$\|R_{\delta}f_{\delta}-f^{'}\|\leq\eta{(\delta)}\longrightarrow{0}$
as $\delta\longrightarrow{0}$.\\  This problem is equivalent to
stably solving the equation
\begin{equation}\label{E3:int}
Au:=\int_{0}^{x}u(t)dt=f(x), \quad
A:H:=L^{2}[0,1]\longrightarrow{L^{2}[0,1]};\quad f(0)=0,
\end{equation}
if noisy data $f_{\delta}$ are given in place of $f$.  In this
case, finding $f^{'}=A^{-1}f$, given the data $f_{\delta}$ is an
ill-posed problem, since equation (\ref{E3:int}) may have no
solution in $L^{2}[0,1]$ if $f_{\delta}\in{L^{2}}[0,1]$ is
arbitrary, subject to only the restriction
$\|f_{\delta}-f\|\leq{\delta},$ and if $f_{\delta}\in{C^{1}[0,1]}$
then $f^{'}_{\delta}$ may differ from $f^{'}$ as much as one
wishes however small $\delta$ is. Also, if $A$ is a linear compact
operator then $A^{-1}$, if it exists is unbounded and hence
equation (\ref{E3:int}) is ill-posed.  The problem is: given
$\{\delta,A,f_{\delta}\}$, find a \emph{stable approximation}
$u_{\delta}$ to the solution $u(x)=f^{'}(x)$ of
 the equation (\ref{E3:int}) in the sense the error estimate
\begin{equation}\label{E4:unois}
\|u_{\delta} - u\| \leq \eta(\delta)\longrightarrow{0} \quad as
\quad \delta\longrightarrow{0}.
\end{equation}
For this we try to construct an operator
$R_{\alpha}:H\longrightarrow{H}$ such that
\begin{equation}\label{E5:reg}
u_{\delta}:=R_{\alpha{({\delta})}}f_{\delta}
\end{equation}
satisfies the error estimate (\ref{E4:unois}).  $R_{\alpha}$
depends on a parameter ${\alpha}$ and is called a
\emph{regularizer} if $R_{\alpha}$ is applicable to any
$f_{\delta}\in{Y}$ and if there is a choice of the regularizing
parameter $\alpha{\equiv\alpha{(\delta})}\rightarrow{0}$ as
$\delta\rightarrow{0}$ such that\\
$R_{\delta}f_{\delta}:=R_{\alpha{({\delta})}}f_{\delta}\longrightarrow{u}\quad{as}\quad{\delta}\longrightarrow{0}$.\\

\emph{ Example 2. The Cauchy problem for the Laplace
equation}\label{Ex2:eg}

We consider the classical problem posed by J.Hadamard.  It is
required to find the solution $u(x,y)$ of the Laplace equation
\begin{equation}\label{E6:cauchy}
u_{xx}+u_{yy}=0
\end{equation}
in the domain $\Omega = \{(x,y)\in R^2 : y>0\}$ satisfying the boundary conditions
\begin{equation}\label{E7:bc}
u(x,0)=0,\quad u_y{(x,0)}=\phi{(x)}=A_{n}\sin nx;\quad
A_{n}\longrightarrow{0}\quad as\quad n\longrightarrow{\infty}.
\end{equation}
The \emph{Cauchy problem} consists of finding a solution of the
equation (\ref{E6:cauchy}) satisfying the conditions
(\ref{E7:bc}). The data differ from zero as little in the sup-norm
as can be wished, if $n$ is sufficiently large.  Its solution is
given by
\begin{equation}\label{E8:soln}
u(x,y)={\frac{A_n}{n}}\sin{(nx)}\sinh{(ny)},
\end{equation}
which, if $A_n = 1/n$, is very large for any value of $y>0$,
because $\sinh{(ny)}=0(e^{ny})$.  As $n\longrightarrow{\infty}$,
the Cauchy data tend to zero in $C^{1}(\mathbb{R})$, and $u\equiv
0$ is a solution to equation (\ref{E6:cauchy}) with
$u=u_{y}=0\quad{at}\quad y=0$.
 Thus, even though the solution to the Cauchy problem (\ref{E6:cauchy})-(\ref{E7:bc}) is unique, continuous dependence of
the solution on the data in the sup-norm does not hold.  This shows that the
Cauchy problem for the Laplace equation (\ref{E6:cauchy}) is an ill-posed problem.\\

\emph{Example 3. Fredholm integral equations of the first
kind}\label{Ex3:eg}

Consider the problem of finding the solution to the integral
equation
\begin{equation}\label{E9:fred}
Au(x)=\int_{0}^{1}K(x,y)u(y)dy=f(x),\quad  0\leq{x}\leq{1}
\end{equation}
where the operator $A:H:=L^{2}(0,1)\longrightarrow{L^{2}(0,1)}$ is
compact and $>0$ almost everywhere, with kernel $K(x,y)$
satisfying the condition:
\begin{equation}\label{E10:frd}
\int_{0}^{1}{\int_{0}^{1}{\left|{K(x,y)}\right|}^{2}}dxdy <
\infty.
\end{equation}

Then $A:H\longrightarrow{H}$ is compact.  A compact operator in an
infinite-dimensional space cannot have a bounded inverse.  That is
the problem (\ref{E9:fred}) is ill-posed.
\newpage
\section{Regularizing family}
Consider the operator equation given by (\ref{E1:oper}) with the
following assumptions:
\begin{enumerate}
\item $A$ is not continuously invertible.
\item For exact values of $f$ and $A$, there exists a solution $u$ of equation (\ref{E1:oper}).
\item $A$ is known exactly, and instead of $f$, we are given its approximation $f_{\delta} \in Y$ such that the estimate (\ref{E2:fnois}) is
satisfied in $Y$.
\end{enumerate}
where ${\delta}>0$ is a numeric parameter characterizing
the errors of input data $f_{\delta}$.

We need a numerical algorithm for solving the operator equation
satisfying the condition that the smaller the value of ${\delta}$
, the closer the approximation to $u$ is obtained, i.e., the error
estimate (\ref{E4:unois}) is satisfied.  The \emph{Regularizing
Algorithm (RA)} is the operator $R_{\alpha}:Y \longrightarrow{X}$
which, for a suitable choice of
${\alpha}\equiv\alpha{({\delta})}$, puts into correspondence to
any pair  $\{\delta, f_{\delta}\}$, the element $u_{\delta}\in X$
such that the error estimate (\ref{E4:unois}) is satisfied where
$u_{\delta}:=R_{\alpha{({\delta})}}f_{\delta}$. For a given set of
data, $u_{\delta}$ is the approximate solution of the problem.
Based on the existence and construction of $RA$, all ill-posed
problems may be classified into regularizable (i.e., the ones for
which a RA exists) and non-regularizable, and \emph{solving an
ill-posed problem means constructing RA for such a problem}.

Let $Au=f$, where $A:X\longrightarrow{Y}$ is a linear injective
operator, $R(A)\neq{\overline{R(A)}}, f \in{R(A)}$ is not known
and the data are the elements $f_{\delta}$ such that the estimate
(\ref{E2:fnois}) is satisfied.  The objective is to find stable
approximation $u_{\delta}$ to the solution $u$ such that the error
estimate (\ref{E4:unois}) is satisfied.  Such a sequence
$u_{\delta}$ is called a \emph{stable solution} to the equation
(\ref{E1:oper}) with the perturbed (or noisy) data.

Let the operator equation $Au=f$, be given, and $f_{\delta}$,
$\|f_{\delta}-f\|\leq{\delta}$, the noisy data be given in place
of $f$.  Let $A$ be injective.  Then $R_{\alpha}$ is called a
\emph{regularizer} of the operator equation if $D(R_{\alpha})=H$
and  there exists
$\alpha{({\delta})}\rightarrow{0},\quad\delta\rightarrow{0}$ such
that $\|R_{\alpha{({\delta})}}f_{\delta}-u\|\rightarrow{0}$ as
$\delta\rightarrow{0}$ for all $u\in{H}$.
\begin{lemma}\label{l33:rh}
If there exists an operator $R_{\alpha}$, $D(R_{\alpha})=H$,
$\|R_{\alpha}\|\leq{a({\alpha})}$, such that
$\|R_{\alpha}Au-u\|:=\eta{({\alpha})}\rightarrow{0}$ as
$\alpha\rightarrow{0}$, and the function
$g({\alpha}):=-\frac{\eta^{'}({\alpha})}{a^{'}{({\alpha})}}$ is
monotone for $\alpha\in{(0,\alpha_{0})},$ $\alpha_{0}>0$ a small
number with $g(+0)=0$, then there exists an
$\alpha{\equiv\alpha}({\delta})\rightarrow{0}$ as
$\delta\rightarrow{0}$ such that
$\|R_{\alpha({\delta})}f_{\delta}-u\|\rightarrow{0}$ as
$\delta\rightarrow{0}$ for all $u\in{H}$.
\end{lemma}
\emph{Proof}: Consider, $\|R_{\alpha}f_{\delta}-u\|$
$\leq{\|R_{\alpha}(f_{\delta}-f)\|}+\|R_{\alpha}f-u\|$
$\leq\delta{a({\alpha})}+\eta{({\alpha})}$ and
$\alpha{({\delta})}$ can be chosen suitably.  One can choose
$\alpha{({\delta})}$ so that the minimization problem

$\delta{a({\alpha})}+\eta{({\alpha})}=min.\longrightarrow{0}$ as
$\alpha\longrightarrow{0}$. \qquad (*)

Equation, $\delta{a^{'}({\alpha})}+\eta^{'}{({\alpha})}=0$. \qquad
(**) is a necessary condition for min. in (*).  Since problem
(\ref{E1:oper}) is ill-posed, one has
$a({\alpha})\rightarrow{\infty}$ as $\alpha\rightarrow{0}$. The
function $a{({\alpha})}$ can be assumed monotone decreasing and
$\eta{({\alpha})}$ can be assumed monotone increasing on
$(0,\alpha_{0})$, $\eta{({0})}=0$.  We assume in lemma,
$\eta{({\alpha})}\rightarrow{0}$, as $\alpha\rightarrow{0}$
[$\eta{({\alpha})}>0$, ,$a^{'}{({\alpha})}<0,
\eta^{'}{({\alpha})}>0,\alpha{>0}$]. Since
$g({\alpha}):=-\frac{\eta^{'}({\alpha})}{a^{'}{({\alpha})}}$ is a
monotone function for $\alpha\in{(0,\alpha_{0})},$ equation (*)
has a unique solution $\alpha\equiv{\alpha{({\delta})}}$ for any
sufficiently small $\delta{>0}$.

More precisely for any fixed $f_{\delta}$, an operator
$R_{\delta}:H\longrightarrow{H}$ such that $R_{\delta}:=
R_{\alpha{({\delta})},{\delta}}$, which depends on ${\delta}$ is a
regularizer for equation (\ref{E1:oper}), if for some choice of
${\alpha},\quad{\alpha}\equiv{{\alpha}({\delta})}$, one has,
\begin{equation}\label{E16:reg}
\|R_{\delta}f_{\delta}-A^{-1}f\|\longrightarrow{0}\quad as \quad
\delta{\longrightarrow{0}},\quad{for}\quad{all}\quad{f}\in{R(A)}.\quad
\end{equation}
so that $u$ is stably approximated.  $\Box$

If such a family $R_{\delta}$ is known then the function
\begin{equation}\label{E17:aprx}
u_{\delta}:=R_{\delta}f_{\delta},
\end{equation}
satisfies the error estimate (\ref{E4:unois}) in view of
(\ref{E16:reg}), i.e., formula (\ref{E17:aprx}) gives a stable
approximation to the solution $u$ of equation (\ref{E1:oper}). The
scalar parameter $\alpha$ is called the \emph{regularization
parameter}.

A construction of a family $R_{{\alpha},\delta}$ of operators such
that there exists a unique $\alpha{(\delta)}$ satisfying
(\ref{E16:reg}) is always used for solving an ill-posed problem
(\ref{E1:oper}).
 The operator $A^{-1}$ is said to be \emph{regularizable}, if there exists a
regularizer $R_{\delta}$ which approximates $A^{-1}$ in the sense
of (\ref{E16:reg}) using the noisy data $\{{\delta},
f_{\delta}\}$.  In the case of well-posed problems, $A^{-1}$ is
always regularizable: one may take $R_{\delta}=A^{-1}$ for all
$\delta$.  This can happen only for well-posed problem.  If the
problem is ill-posed then there does not exist a regularizer
independent of the noise $\delta$.

\emph{Example}:  In the stable numerical differentiation example,
we shall take
\begin{equation}\label{E18:diff}
Au(x)= \int_{0}^{x}u(s)ds=f(x), \quad{0}\leq{x}\leq{1}.
\end{equation}
Suppose $f_{\delta}(x)$ is given in place of $f$:
\begin{equation}\label{E19:deri}
\|f-f_{\delta}\|_{L^{\infty}{(0,1)}}\leq{\delta}.
\end{equation}
Following A.G.Ramm \cite{aR68}, we choose a regularizer
$R_{\delta}$ of the form:
\begin{equation}\label{E20:dereg}
u_{\delta}=R_{\delta}f_{\delta}:=\frac{f_{\delta}(x+h({\delta}))-f_{\delta}(x-h({\delta}))}{2h},\quad
h({\delta})=\sqrt{\frac{2\delta}{M}}>0.
\end{equation}

We note that  A.G.Ramm \cite{aR03} has given a new notion of
regularizer.  According to \cite{aR03} a family of operators
$R(\delta)$ is a \emph{regularizer} if\\
$sup\|R(\delta)f_{\delta}-v\|\leq\eta{(\delta)}\longrightarrow{0}
\quad as \quad \delta{\longrightarrow{0}},$\\  where the supremum
is taken over all
$v\in{S_{\delta}=\{v:\|Av-f_{\delta}\|\leq{\delta{}}{,v\in
{K}}}{\}}$, and $K$ is a compactum in $X$ to which the solution
$u$ belongs.  The difference between Ramm's  definition and the
original definition is that in the original definition $u$ is
fixed, one does not know the solution $u\in{K},$ the only
information available is a family $f_{\delta}$ and some apriori
information about the solution $u$, while in the new definition
$v$ is any element of $S_{\delta}$ and the supremum, over all such
$v$, of the norm above must go to zero as $\delta$ goes to zero.
 This definition is more natural in the sense that not only the
solution $u$ to (\ref{E1:oper}) satisfies the estimate
$\|Au-f_{\delta}\|\leq{\delta}$, but many $v\in{K}$ satisfy such
an inequality $\|Av-f_{\delta}\|\leq{\delta},v\in{K}$, and the
data $f_{\delta}$ may correspond, to any $v\in{S_{\delta}}$, and
not only to the solution of problem (\ref{E1:oper}).
\chapter{Review of the methods for solving ill-posed problems}
In this chapter, we shall discuss four different methods,
(variational regularization method, quasi-solutions method,
iterative regularization method and dynamical systems method) for
constructing regularizing families for ill-posed problems
(\ref{E1:oper}) with bounded operators.  See also
A.G.Ramm\cite{aR02} for ill-posed problems with unbounded operators.

\section{Variational regularization method}
This method consists of solving a variational problem, which was
proposed by D.Phillips (1962) and studied by A.N.Tikhonov (1963)
et al by  constructing regularizers for solving ill-posed
problems.

Consider equation (\ref{E1:oper}), which has to be solved, where
$A:X\longrightarrow{Y}$ is assumed to be a bounded, linear,
injective operator, with $A^{-1}$ not continuous, $f \in R(A)$ is
not known, and the data are the elements $\{{\delta}, A,
f_{\delta}\}$, where the noise level ${\delta}>0$ is given,
estimate (\ref{E2:fnois}) holds, and the noisy data $f_{\delta}$
is the ${\delta}-$approximation of $f$.  The problem is: given
$\{{\delta}, A, f_{\delta}\}$, find the stable solution
$u_{\delta}$ such that the error estimate (\ref{E4:unois}) holds.

Let equation (\ref{E1:oper}) have a minimal-norm solution $y$.
Variational regularization method consists of solving the
variational problem (minimization problem) and constructing a
stable approximation to solution $y$ with minimal-norm such that
$y\bot{N(A)}$.  Assume without loss of generality $\|A\|\leq{1}$,
and then $\|A^{*}\|\leq{1}$.  Let $B:=A^{*}A$, then $B\geq{0}$ and
is a bounded, self-adjoint operator. The equation $Au=f$ is
equivalent to the equation $Bu=q$, where $q:=A^{*}f$.  Assume that
$q_{\delta}$ is given in place of $q$,
$\|q-q_{\delta}\|\leq\|A^{*}\|\delta\leq\delta$.  Since
$N(A)=N(B)$, $y\perp{N(B)}$ and
$\|(B+{\alpha})^{-1}\|\leq{\frac{1}{\alpha}}$.

Consider the problem of finding the minimum of the functional
\begin{equation}\label{E21:funcal}
F(u):={\|Au-f_{\delta}\|}^{2}+{\alpha}{\|u\|}^{2}=min.
\end{equation}
where ${\alpha}>0$ is the regularization parameter.  The
functional $F(u)$ is a function of two parameters $\alpha$ and
$\delta$.
 Solutions of variational problem (\ref{E21:funcal}) are called \emph{minimizers}.
First, we shall prove the following two lemmas.
\begin{lemma}\label{l1:min}
\emph{Existence of minimizers:} For arbitrary $\alpha{>}0$ and
$\delta{>}0$, there exists a solution $u_{\alpha,\delta}$ to
variational problem (\ref{E21:funcal}), in the sense
$F(u_{\alpha,\delta})\leq{F(u)}$ for all $u\in{X}$.
\end{lemma}
\begin{lemma}\label{l2:min}
\emph{Uniqueness of minimizers:} The solution $u_{\alpha,\delta}$
of variational problem (\ref{E21:funcal}) is unique.
\end{lemma}
\emph{Proof of Lemma \ref{l1:min}}: Define\\
$m:=\inf_{u\in{H}}{F(u)}$.  Note that $m\equiv{m}(\delta
)\geq{0}$, and $m\leq{F(u)}, u\in{H}$.\\  Let $\{u_{n}\}\in{D(F)}$
be a minimizing sequence for the functional $F$, such that
\begin{equation}\label{E22:min}
m\leq{F(u_{n})}\leq{m+\epsilon_{n}},\quad\epsilon_{n}\longrightarrow{0}\quad{as}\quad{n}\longrightarrow{\infty}.
\end{equation}
So, $\alpha{\|u_{n}\|}^{2}\leq{m+\epsilon_{n}}$.  Hence
 ${\|u_{n}\|}^{2}\leq{\frac{m+\epsilon_n}{\alpha}.}$

Since a bounded set in $H$ contains a weakly convergent
subsequence, there exists a weakly convergent subsequence of
$\{u_n\}$, denoted again by $\{u_n\}$, with
$u_{n}\rightharpoonup{u}$.  This implies by continuity of $A$,
$Au_{n}\rightharpoonup{Au}$.  Thus,

$\|u\|\leq{\liminf_{n\rightarrow{\infty}}{\|u_{n}\|}}$ and
${\|Au-f_{\delta}\|}^{2}\leq{\liminf_{n\rightarrow{\infty}}{\|Au_{n}-f_{\delta}\|}^{2}}$.

So from equation (\ref{E22:min}) we obtain\\
$m\leq{{\|Au-f_{\delta}\|}^{2}+\alpha{{\|u\|}^{2}}}$
$\leq\liminf_{n\rightarrow
\infty}({{\|Au_{n}-f_{\delta}\|}^{2}+{\|u_{n}\|}^{2}})$
$\leq\liminf_{n\rightarrow{\infty}}F(u_{n})$\\
$\leq\lim_{n\to{\infty}}(m+\epsilon_{n})=m$.

So we have $m\leq{F(u)}\leq{m}.$  Hence $F(u)=m$. So $u$ is the
minimizer of $F(u)$.  Thus we have proved the existence of the
minimizer for the variational problem (\ref{E21:funcal}).  $\Box$

\emph{Proof of Lemma \ref{l2:min}}: Since $u$ is the minimizer of
variational problem (\ref{E21:funcal}), it follows that
$F(u)\leq{F(u+\epsilon\eta)}$, for any $\eta\in{H}$, and for any
$\epsilon\in{(0,\epsilon_{0})}$. So,
$\lim_{\epsilon\to{0}}{\frac{F(u+\epsilon\eta)-F(u)}{\epsilon}}\geq{0}$.

Assuming that $F^{'}(u)$ exists, this implies that
$(F^{'}(u),\eta)\geq{0}$ for all $\eta\in{}H$ and hence
$F^{'}(u)=0$.  We shall calculate the derivative of
(\ref{E21:funcal}) with respect to  $\epsilon$ at $\epsilon=0$ and
get:

$\frac{d}{d\epsilon}F(u+\epsilon\eta)|_{\epsilon=0}$
$=[\frac{d}{d\epsilon}{{\|Au+A\epsilon\eta-f_{\delta}\|}^{2}+\alpha{\|u+\epsilon\eta\|}^{2}}]|_{\epsilon=0}$.

Since

$[\frac{d}{d\epsilon}{\|u+\epsilon\eta\|}^{2}]|_{\epsilon=0}$
$=[\frac{d}{d\epsilon}(u+\epsilon\eta,
u+\epsilon\eta)]|_{\epsilon=0}$
$=[(\eta{,}u+\epsilon\eta),(u+\epsilon\eta,\eta)]|_{\epsilon=0}$
$=(\eta{,}u)+(u,\eta)$ $=2Re(u,\eta)$

and

$[\frac{d}{d\epsilon}{\|Au+A\epsilon\eta-f_{\delta}\|}^{2}]|_{\epsilon=0}$
$=[\frac{d}{d\epsilon}(Au+\epsilon{A}\eta{-}f_{\delta},Au+\epsilon{A}\eta{-}f_{\delta})]|\epsilon_{=0}$
$=(A\eta{,}Au-f_{\delta})+(Au-f_{\delta},A\eta{)}$
$=2Re(A^{*}Au-A^{*}f_{\delta},\eta).$\\
So,\\
$\frac{d}{d\epsilon}F(u+\epsilon\eta)_{|\epsilon=0}$
$=2Re(A^{*}Au-A^{*}f_{\delta}+\alpha{u},\eta)=0$\\   for all
$\eta\in{H}$.  Hence we obtain,
\begin{equation}\label{E23:con}
A^{*}Au+\alpha{u}=A^{*}f_{\delta}.
\end{equation}
Thus if $u$ is a minimizer of $F(u)$, then equation
(\ref{E23:con}) holds.  We claim that equation (\ref{E23:con}) has
not more than one solution.  For this, it is sufficient to prove
that $A^{*}Aw+\alpha{w}=0$ implies $w=0$.  Suppose that
\begin{equation}\label{E24:homo}
A^{*}Aw+\alpha{w}=0,\quad\alpha{=}constant>0.
\end{equation}
then
$0=((A^{*}A+\alpha)w,w)$\\
$=(A^{*}Aw,w)+\alpha{(w,w)}$ $=(Aw,Aw)+\alpha{(w,w)}$
$={\|Aw\|}^{2}+\alpha{\|w\|}^{2}\geq\alpha{\|w\|}^{2},\quad
\alpha{>0}.$

Therefore, $w=0$. Hence the solution to equation (\ref{E23:con})
is unique, and is given by the formula
\begin{equation}\label{E25:sol}
u_{\alpha{,}\delta}:= (A^{*}A+\alpha)^{-1}A^{*}f_{\delta}.
\end{equation}
for every $\alpha>0$ and the operator $(A^{*}A+\alpha)^{-1}$
exists and is bounded by
$\|(A^{*}A+\alpha)^{-1}\|\leq\frac{1}{\alpha}$, because
$A^{*}A\geq{0}$.  $\Box$

We shall now consider the main theorem which gives us a method for
constructing a regularizing family for the ill-posed problem
(\ref{E1:oper}).
\begin{theorem}\label{t1:main}
Assume that $A$ is a linear bounded, injective operator, equation
$Au=f$ is solvable, $A^{-1}$ is not continuous, and let $y$ be the
minimal-norm solution: $Ay=f$, $y\bot{N(A)}$.  If $f_{\delta}$ is
given in place of $f$, $\|f-f_{\delta}\|\leq{\delta}$.  Then for
any $\alpha{>0}$, minimization problem (\ref{E21:funcal}) has a
unique solution $u_{\alpha{,\delta}}$ given by the formula
(\ref{E25:sol}). Moreover, if $\alpha\equiv\alpha{(\delta)}$ is
such that
\begin{equation}\label{E26:condi}
\alpha{(\delta)}\longrightarrow{0} \quad and \quad
\frac{\delta}{\alpha{(\delta)}}\longrightarrow{0} \quad as \quad
\delta\longrightarrow{0},
\end{equation}
then $u_{\delta}:=u_{\alpha{(\delta}){,\delta}}\longrightarrow{y}$
\quad as \quad $\delta\longrightarrow{0}$.
\end{theorem}
\emph{Proof:} The proofs of existence and uniqueness of the
minimizers were given above.  Let us prove the last conclusion of
the theorem. Assume that condition (\ref{E26:condi}) holds.

Define the regularizer (by means of formula (\ref{E25:sol}) so
that it satisfies equation (\ref{E21:funcal})):
\begin{equation}\label{E27:reg}
R_{\alpha}f_{\delta}:=R_{\alpha{({\delta})}}f_{\delta}:=u_{{\alpha},{\delta}}=(B+{\alpha})^{-1}A^{*}f_{\delta},\quad\alpha{=\alpha{({\delta})}}
.
\end{equation}
If $f=Ay$, then
$R{_\alpha}Ay=(B+{\alpha})^{-1}By\longrightarrow{y}$ as
$\alpha\longrightarrow{0}$.  Here the assumption $y\perp{N(A)}$ is
used.  We claim that
\begin{equation}\label{E28:nor}
\|R_{\alpha}f_{\delta}-y\|\longrightarrow{0}\quad{as}\quad
\alpha\longrightarrow{0}.
\end{equation}
Since,

$\|R_{\alpha}f_{\delta}-y\|\leq\|R_{\alpha}(f_{\delta}-f)\|+\|R_{\alpha}f-y\|$
$\leq{\|f_{\delta}-f\|}\|R_{\alpha}\|+\eta{(\alpha)}$
$\leq{\frac{\delta}{\alpha}}\|A^{*}\|+\eta{({\alpha})}$
$\leq{\frac{\delta}{\alpha}}+\eta{({\alpha})}\rightarrow{0}$ as
$\delta\longrightarrow{0}$, because of (\ref{E26:condi}).  $\Box$

\newpage
\section{Discrepancy principle for variational regularization method}

\emph{Assumptions (A)}   Let $A$ be a linear, bounded, injective
operator.  Let equation $Au=f$ be solvable.  Let $A^{-1}$ not be
continuous and let there exist a minimal-norm solution $y$ such
that $Ay=f$, $y\perp{N(A)}$. Let $f_{\delta}$ be given in place of
$f$,  $\|f-f_{\delta}\|\leq\delta$.  Let $u_{\delta}$ be the
stable solution of equation (\ref{E1:oper}) given by
$u_{\delta}={(B+{\alpha})}^{-1}A^{*}f_{\delta}$.

The discrepancy principle (DP) introduced by Morozov is used as an
aposteriori choice of the regularization parameter $\alpha$ and
this choice yields convergence of the variational regularization
method.  Choose the regularization parameter
$\alpha{=}\alpha{({\delta})}$ as the root of the equation

$\|Au_{\alpha{,\delta}}-f_{\delta}\|=C\delta{,}
\quad{C}=constant>1.$   \qquad  (*)\\
The above equation, is a non-linear equation with respect to
${\alpha}$.  It defines $\alpha$ as an implicit function of
$\delta$.  Let us assume that $\|f_{\delta}\|>C\delta.$

\begin{theorem}\label{td:dp}
Suppose that the assumptions (A) above holds.  Then there exists a
unique solution $\alpha{=\alpha{(\delta})}>0$ to equation (*) and
$\lim_{{\delta}\to{0}}\alpha{(\delta})=0$.
 Moreover, if
$u_{\delta}:=u_{\alpha{(\delta}),{\delta}}$ is given by formula
(\ref{E27:reg}), then
$\lim_{\delta\rightarrow{0}}\|u_{\delta}-y\|=0.$
\end{theorem}
\emph{Proof:}  Denote $Q:=AA^{*}.$  Then, $N(Q)=N(A^{*})$, and
${(B+{\alpha})}^{-1}A^{*}=A^{*}{(Q+{\alpha})}^{-1}$.  By
variational regularization method, for any $\alpha{>0}$,
minimization problem (\ref{E21:funcal}) has a unique solution
$u_{{\alpha},\delta}={(B+{\alpha})}^{-1}A^{*}f_{\delta}$.  Since,
$\|f_{\delta}-f\|\leq{\delta}$; $f=Ay$, so, $f\perp{N(A^{*})}$.
From equation (*),

$C^{2}{\delta}^{2}$
$=\|Au_{{\alpha},\delta}-f_{\delta}\|^{2}$
$=\|[A(B+{\alpha})^{-1}A^{*}-I]f_{\delta}\|^{2}$
$=\|[Q(Q+{\alpha})^{-1}-I]f_{\delta}\|^{2}$

$={\alpha}^{2}\|(Q+{\alpha})^{-1}f_{\delta}\|^{2}$
$={\alpha}^{2}\int_{0}^{\|Q\|}\frac{d(E_{\lambda}f_{\delta},f_{\delta})}{(\lambda{+\alpha})^{2}}:=I(\alpha{,\delta})$

one has,\\
$\lim_{\alpha\rightarrow{\infty}}I(\alpha{,\delta})=\int_{0}^{\|Q\|}d(E_{\lambda}f_{\delta},f_{\delta})=\|f_{\delta}\|^{2}>C^{2}{\delta}^{2}$.

and\\
$\lim_{\alpha\rightarrow{+0}}I(\alpha{,\delta})=\|P_{N(Q)}f_{\delta}\|^{2}$.

where $P_{N(Q)}$ is the orthogonal projection onto the null-space
of $Q$,\\
$P_{(a,b)}f_{\delta}:=\int_{a}^{b}dE_{\lambda}f_{\delta}$\\
 and\\
$\|P_{(a,b)}f_{\delta}\|^{2}:=\int_{a}^{b}d(E_{\lambda}f_{\delta},f_{\delta})$.\\
So\\
$\|P_{N(Q)}f_{\delta}\|^{2}\leq\int_{0}^{\epsilon}\frac{{\alpha}^{2}{d(E_{\lambda}f_{\delta},f_{\delta})}}{({\alpha{+\lambda}})^{2}}$
$\leq\int_{0}^{\epsilon}d(E_{\lambda}f_{\delta},f_{\delta})=\|P_{(0,{\epsilon})}f_{\delta}\|^{2}$.\\
Since\\
$P_{N(A^{*})}f_{\delta}=P_{N(A^{*})}(f_{\delta}-f)+P_{N(A^{*})}f$\\
and since $P_{N(A^{*})}f=0$, one has

$\|P_{N(A^{*})}f_{\delta}\|^{2}=\|P_{N(A^{*})}(f_{\delta}-f)\|^{2}$
$\leq\|f_{\delta}-f\|^{2}\leq{\delta}^{2}$.

Thus,

$\lim_{\alpha\rightarrow{+0}}I(\alpha{,\delta})=\|P_{N(Q)}f_{\delta}\|^{2}$
$=\|P_{N(A^{*})}f_{\delta}\|^{2}\leq{\delta}^{2}<C^{2}{\delta}^{2}$.

Equation (*) is a non-linear equation of the form
$C^{2}{\delta}^{2}=I{({\alpha},{\delta})}$, for a given fixed pair
$\{f_{\delta},\delta\}$, the function $I{({\alpha},{\delta})}$
satisfies
$\lim_{\alpha\rightarrow{+0}}I(\alpha{,\delta})<{C^{2}}{\delta}^{2}$
and
$\lim_{\alpha\rightarrow{\infty}}I(\alpha{,\delta})>{C^{2}}{\delta}^{2}$.
Hence $I{({\alpha},{\delta})}$ is a monotone increasing function
of $\alpha$ on $(0,{\infty})$.  Hence,  equation (*) has a unique
solution $\alpha{=}\alpha{({\delta})}$.  Now let us prove that
$\lim_{{\delta}\rightarrow{0}}{\alpha}({\delta})=0$.  Suppose that
$\alpha{({\delta})}\geq{\alpha_{0}}>0$.  So as
$\delta{\rightarrow{0}}$,

$0\leftarrow{C^{2}{\delta}^{2}}=\int_{0}^{\|Q\|}\frac{{\alpha}^{2}{d(E_{\lambda}f_{\delta},f_{\delta})}}{({\alpha{+\lambda}})^{2}}$
$\geq\int_{0}^{\|Q\|}\frac{{\alpha_{0}}^{2}{d(E_{\lambda}f_{\delta},f_{\delta})}}{({\alpha_{0}{+\lambda}})^{2}}$
$\geq\frac{{\alpha_{0}}^{2}}{({\alpha_{0}{+\|Q\|}})^{2}}\|f_{\delta}\|^{2}>0$.

This contradicts the assumption that $\alpha{({\delta})}>0$.  It
remains to prove the last conclusion of the theorem.
 Define the regularizer by the formula (by means of formula
 (\ref{E27:reg}))\\
$R_{\delta}f_{\delta}:=u_{\delta}:=u_{\alpha{(\delta}),\delta}$,\\
where $\alpha{(\delta})$ is given by the discrepancy principle.
Let us prove that $\|u_{\delta}-y\|\rightarrow{0}$ as
$\delta\rightarrow{0}$.  Since $u_{\delta}$ is a minimizer of (\ref{E21:funcal}), we have\\
$F(u_{\delta})\leq{F(u)},\quad$
$\|Au_{\delta}-f_{\delta}\|^{2}+{\alpha{({\delta})}}\|u_{\delta}\|^{2}$
$\leq{\delta}^{2}+\alpha{({\delta})}{\|y\|}^{2}$  \qquad
(**)\\
then from equations (*) and (**), we obtain,\\
$\|u_{\delta}\|^{2}\leq\|y\|^{2},\quad\|u_{\delta}\|\leq\|y\|$.

This implies that there exists $v$ such that
$u_{\delta}\rightharpoonup{v}\quad
as\quad\delta\longrightarrow{0}$ and by continuity of $A$,
$Au_{\delta}\longrightarrow{Av}.$  So from (**), as
${\delta}\rightarrow{0}$ and $\alpha\rightarrow{0}$ it follows
that  $Au_{\delta}\longrightarrow{f}.$  So,  $Av=f$.  Since $A$ is
injective, this implies  that $v=y$.  So,
$u_{\delta}\rightharpoonup{y}\quad as \quad
\delta\longrightarrow{0}$.  Also, since,
$\|u_{\delta}\|\leq\|y\|$,\\
$\|y\|\leq\liminf_{\delta\rightarrow{0}}\|u_{\delta}\|\leq\limsup_{\delta\rightarrow{0}}\|u_{\delta}\|\leq\|y\|$.\\
Therefore, $\lim_{\delta\rightarrow{0}}\|u_{\delta}\|$ exists and
$\lim_{\delta\rightarrow{0}}\|u_{\delta}\|=\|y\|$.  Thus,
$\lim_{\delta\rightarrow{0}}\|u_{\delta}-y\|=0$.  Hence the
theorem is proved.  $\Box$

Note that: $F(P_{N(A)^{\perp}}u)\leq{F(u)}$.

Proof: Let $u=u_{0}+u_{1}$, where $u_{0}\in{N(A)}$ and
$u_{1}\in{N(A)^{\perp}}$.

$Au=Au_{0}+Au_{1}=Au_{1}$,   since $Au_{0}=0$.

So, $\|Au-f_{\delta}\|^{2}=\|Au_{1}-f_{\delta}\|^{2}$.\\ Also,
$\alpha\|u\|^{2}=\alpha{[\|u_{0}\|^{2}+\|u_{1}\|^{2}]}\geq\alpha\|u_{1}\|^{2}$.

This implies that,  $F(u_{1})\leq{F(u)}$. So
$F(u_{1})=F(P_{N(A)^{\perp}}u)\leq{F(u)}$.  Hence a minimizer of
$F$ is necessarily orthogonal to null-space of $A$.

\emph{Remark:} A.G.Ramm \cite{aR03a1} has generalized the
discrepancy principle for the cases: (a) when $A$ is not
injective, (b) when $A$ is not compact and not injective and (c)
when $A^{-1}$ is not continuous.  He has also shown that
discrepancy principle, in general does not yield convergence which
is uniform with respect to the data.
\newpage
\section{The method of quasi-solution}
The method of quasi-solution was given by Ivanov (1962).  It is
similar to the variational regularization method except that there
is a restriction on the functional defined.

Consider the operator equation (\ref{E1:oper}) which has to be
solved, where $A$ is assumed to be a bounded, linear injective
operator on Banach spaces $X$ and $Y$ or $R(A)$ is assumed not to
be closed, so that the problem is ill-posed.   The data are the
elements $\{{\delta}, A, f_{\delta}\}$, where the noise level
${\delta}>0$ is given such that the estimate (\ref{E2:fnois})
holds, i.e., the noisy data $f_{\delta}$ is the
${\delta}-$approximation of $f$. The problem is: given
$\{{\delta}, A, f_{\delta}\}$, find the stable solution
$u_{\delta}$ such that the error estimate (\ref{E4:unois}) holds.
Let equation (\ref{E1:oper}) have a solution $y\in{K}$, a convex
compactum (closed, pre-compact subset) of $X$.  Consider the
variational problem:
\begin{equation}\label{E30:funcal}
F(u):=\|{Au-f_{\delta}\|}\longrightarrow{inf.}, \quad{u}\in{K}.
\end{equation}
\begin{definition}\label{D2:quasi}
A quasi-solution of equation (\ref{E1:oper}) on a compactum $K$ is
a solution to the minimization problem (\ref{E30:funcal}).
\end{definition}
\begin{lemma}: Existence of quasi-solution:\\
Assume that $A$ is a bounded linear injective operator and that
equation (\ref{E1:oper}) holds. Assume that equation
(\ref{E1:oper}) has a solution $y\in{K}$ a compactum of X.  Then
the minimization problem (\ref{E30:funcal}) has a stable solution
$u_{\delta}\in{K}$ such that
$\|u_{\delta}-y\|\longrightarrow{0}\quad{as}\quad\delta\longrightarrow{0}$.
\end{lemma}
\emph{Proof}: Denote
\begin{equation}\label{E31:min}
m({\delta}):=\inf_{u\in{K}}{\|Au-f_{\delta}\|}.
\end{equation}
Since, the infimum $m=m({\delta})$ depends on $f_{\delta}$ and
since $y\in{K},$ we have,

$m({\delta})=\inf_{u\in{K}}{\|Au-f_{\delta}\|{\leq}\|Ay-f_{\delta}\|}$
$=\|f-f_{\delta}\|{\leq}{\delta}$.\\  So
$m({\delta})\longrightarrow{0}\quad{as}\quad
\delta{\longrightarrow{0}}.$  Let ${u_{n}}$ be a minimizing
sequence in $K$:
\begin{equation}\label{E32:min}
F(u_{n}):=\|{Au_{n}-f_{\delta}\|}\longrightarrow{m({\delta})},
\quad u_{n}\in{K},  \quad n\longrightarrow{\infty}.
\end{equation}
So we have  $\sup{_n}\|Au_{n}\|<{\infty}.$  Let us now take
$\delta\longrightarrow{0}$.  Since ${u_{n}}\in{K}$ and $K$ is a
compactum, there exists a convergent subsequence in $K$, which we
again denote by $u_{n}$, such that
$u_{n}\longrightarrow{u_{\infty}}$.   Since $K$ is a compactum, it
is closed.  Therefore the limit $u_{\infty}\in{K}$.   By
continuity of $A,$ this implies that $
Au_{n}\longrightarrow{Au_{\infty}},$ and
$\lim_{n\to{\infty}}\|Au_{n}-f_{\delta}\|=\|Au_{\infty}-f_{\delta}\|=m(\delta{)}.$

Denote
\begin{equation}\label{E33:quasi}
u_{\delta}:=u_{\infty}\equiv{u_{\infty}}{(\delta)}\in{K}.
\end{equation}
Thus, $u_{\delta}$ is the solution of the minimization problem
(\ref{E30:funcal}):

$\|Au_{\delta}-f_{\delta}\|=m(\delta{)}$.

It remains to be shown that $u_{\delta}\in{K}$, is the
quasi-solution of equation (\ref{E1:oper}).  Now, as
${\delta}\longrightarrow{0}$, there exists a subsequence
$u_{\delta_{n}}\in{K}$ which is again denoted by $u_{n},$ such
that $u_{n}\longrightarrow{v}\in{K}$.  By continuity of $A$, this
implies that $Au_{n}\longrightarrow{Av}.$  Therefore, since
$m({\delta})\longrightarrow{0}\quad{as}\quad\delta\longrightarrow{0},$
\begin{equation}\label{E34:quasi}
\|Av-f\|\longrightarrow{0},\quad
as\quad\delta{\longrightarrow{0}}.
\end{equation}
Since, $A$ is injective, $v=y$.  Thus,
\begin{equation}\label{E35:quas}
\lim_{{\delta}\to 0}\|{u_{n}}-y\|{=0}.
\end{equation}
Since the limit y of any subsequence $u_{n}$ is the same, the
whole sequence  $u_{n}$ converges to $y$.  Thus a quasi-solution
exists.  Hence lemma 5 is proved.  $\Box$

It remains to be proved the uniqueness and its continuous
dependence on $f$ of the quasi-solution.
\begin{theorem}\label{t2:main}
If $A$ is linear, bounded and injective operator, $K$ is a convex
compactum and the functional $F(u)$ in minimization problem
(\ref{E30:funcal}) is strictly convex, then for any $f$, the
quasi-solution exists, is unique, and depends on $f$ continuously.
\end{theorem}
\emph{Proof}:\\ The following lemmas, are needed for the proof of
theorem (3).
\begin{lemma}
Let $\inf_{u\in{K}}\|u-f\|:=dist(f,K):=m(f)$.  Then there exists a
unique element $u\equiv{u(f)}:=P_{K}f\in{K}$ called the metric
projection of $f$ onto $K$ such that $\|P_{K}f-f\|=dist(f,K)$.
\end{lemma}
\emph{Proof}: \emph{Existence of $P_{K}f$}:  Let ${u_{n}}$ be a
minimizing sequence in $K$,  $\|u_{n}-f\|\longrightarrow{m(f)}$,
Let $n\longrightarrow{\infty}$. Then there exists a convergent
subsequence in $K$, which we again denote by $u_{n}$, such that
$u_{n}\longrightarrow{u}\in{K}$.  Thus, $\|u-f\|=m(f)$. So
$u=P_{K}f$.

\emph{Uniqueness of $P_{K}f$}:  Suppose there exists $u,v$ which
are distinct metric projections.  Then
$m(f)=\|u-f\|=\|v-f\|\leq\|w-f\|$ for all $w\in{K}$. Since $K$ is
convex, $\frac{u+v}{2}\in{K}$.  This implies that

$m(f)\leq\|\frac{u+v}{2}-f\|=\|\frac{u-f+v-f}{2}\|\leq\frac{\|u-f\|+\|v-f\|}{2}=m(f).$

So,\\ $\|\frac{u+v}{2}-f\|=m(f)$.

Thus,\\ $\|u-f\|=\|v-f\|=\|\frac{u-f+v-f}{2}\|$.

Since $X$ is strictly convex, it follows that
$(u-f)=\lambda{(v-f)}$. Since $\|u-f\|=\|v-f\|$, $\lambda=+1,-1$.
If $\lambda=1$, then $u=v$ which is a contradiction.  If
$\lambda=-1$, then $f=\frac{u+v}{2}$.  Since $K$ is convex, this
implies that $f\in{K}$, this gives that $P_{K}f=f$ which is a
contradiction.  Thus $P_{K}$ is a bijective mapping onto $K$.
Hence Lemma 6 is proved.  $\Box$
\begin{lemma}
$dist(f,K)$ is a continuous function of $f.$
\end{lemma}
\emph{Proof}:  Let $dist(f,K):=m(f)$.  Suppose $f\rightarrow{g}$.
Then to prove that $m(f)\rightarrow m(g)$.  Let
$u(f)=P_{K}f\in{K}$ and $u(g)=P_{K}g\in{K}$. ($f$ and $g$ are
arbitrary they need not be in $K$).\\ So

$\|u(f)-f\|=inf_{u\in{K}}\|u-f\|$ and $\|u(g)-g\|=inf_{u\in{K}}\|u-g\|$.

So,\\ $m(f)=\|u(f)-f\|\leq\|u(g)-f\|\leq\|u(g)-g\|+\|g-f\|$.

Hence, $m(f)-m(g)\leq\|g-f\|$.\\ Similarly,\\
$m(g)-m(f)\leq\|g-f\|$.\\ Thus,\\ $|m(f)-m(g)|\leq\|f-g\|$.\\
Hence lemma 7 is proved.  $\Box$

\begin{lemma}
$P_{K}f$ is a continuous function of $f$ (in a strictly convex Banach space).
\end{lemma}
\emph{Proof}: Suppose there is a sequence $f_{n}\rightarrow{g}$.
Then to prove that

$\|u(f_{n})-u(g)\|\rightarrow{0}\quad{as}\quad{n}\rightarrow{\infty}$.
\qquad (*)

Suppose (*) is not true, so that there is a sequence $u_{n}$ in
$K$ which does not satisfy (*).  Since $K$ is a compactum, there
exists a subsequence ${u_{n}}_{k}\in{K}$ of $u_{n}$, which is
denoted again by $u_{n}$ such that,\\
$\|u_{n}-u(g)\|\geq\epsilon>0$.\\  Also, since $K$ is closed,
$u_{n}\rightarrow{v}\in{K}$, so that

$\|v-u(g)\|\geq\epsilon{>}0$. \qquad   (**)

$\|u(g)-g\|\leq\|v-g\|$. \qquad  (***)

Now,\\ $\|v-g\|\leq\|v-u_{n}\|+\|u_{n}-f_{n}\|+\|f_{n}-g\|$.

By lemma 7, since $f_{n}\rightarrow{g}$,
$\|u_{n}-f_{n}\|=m(f_{n})\rightarrow{m(g)}$.

Also we have,\\  $\|v-u_{n}\|\rightarrow{0}$ and
$\|f_{n}-g\|\rightarrow{0}$.

Thus,\\ $\|v-g\|\leq{m(g)}=\|u(g)-g\|$.  \qquad (****)

So, by inequalities (***) and (****),\\ $\|u(g)-g\|=\|v-g\|$.\\
This by uniqueness implies that $v=u(g)$. This contradicts
inequality (**).  Hence lemma 8 is proved. $\Box$

\begin{lemma}
If $A$ is a closed (possibly non-linear) injective map over a
compactum $K\subset{X}$ onto $AK$, then $A^{-1}$ is a continuous
map of $AK$ onto $K$.
\end{lemma}
\emph{Proof}: Let $f_{n}=Au_{n}$, where the sequences
$u_{n}\in{K}$, and $f_{n}\in{AK}$.  Assume that
$f_{n}\rightarrow{f}$.  Then to prove that $f\in{AK}$, that is to
prove that there exists a $u\in{K}$ such that
$u_{n}=A^{-1}f_{n}\rightarrow{u}=A^{-1}f$.  Since $K$ is a
compactum, and since $u_{n}\in{K}$, there exists a convergent
subsequence, which is again denoted by $u_{n}\in{K}$ such that
$u_{n}\rightarrow{u}$.  Since $K$ is a compactum, it is closed, so
$u\in{K}$.  Because any convergent subsequence of $u_{n}$
converges to a unique limit $u$, implies that the whole sequence
converges to $u$.  Since $u_{n}\rightarrow{u}$,
$f_{n}=Au_{n}\rightarrow{f}$ and $A$ is closed, therefore, $Au=f$.
Since $A$ is injective, this implies that $u=A^{-1}f$.  Hence
lemma 9 is proved.  $\Box$

\emph{Proof of continuous dependence on f in theorem (3)}:

Existence of quasi-solution is proved in \emph{lemma (5)}.  Since
$K$ is convex and $A$ is linear, so $AK$ is convex.  Since $AK$ is
convex and $F$ is strictly convex, by lemma (6), $P_{AK}f$ exists
and is unique.  By lemma (8), $P_{AK}f$ depends on $f$
continuously.  Let $Au=P_{AK}f$.  Since $A$ is injective
$u=A^{-1}P_{AK}f$ is uniquely defined and by lemma (9), depends
continuously on $f$.  Thus theorem (3) is proved.  $\Box$

\emph{Remark}:

By theorem (3), if $K$ is a convex compactum of $X$ which contains
the solution $u$ to equation (\ref{E1:oper}), if $A$ is an
injective linear bounded operator, and $F$ is strictly convex,
then $u_{\delta}=A^{-1}P_{AK}f_{\delta}$ satisfies
$\|u_{\delta}-u\|\rightarrow{0}\quad{as}\quad{\delta\rightarrow{0}}$.
The function $u_{\delta}$ can be found as the unique solution to
the minimization problem (\ref{E30:funcal}) with $f_{\delta}$ in
place of $f$.  Further instead of assuming operator $A$  to be
bounded, $A$ can be assumed to be closed, since a bounded operator
defined everywhere is closed.

\newpage
\section{Iterative regularization method}
Consider the operator equation (\ref{E1:oper}) which has to be
solved, where $A:H\longrightarrow{H}$ is assumed to be a bounded,
linear injective operator on a Hilbert space H with $A^{-1}$
unbounded.  So the problem is ill-posed.  The data are the
elements $\{{\delta}, A, f_{\delta}\}$, where the noise level
${\delta}>0$ is given such that estimate (\ref{E2:fnois}) holds,
i.e., $f_{\delta}$ is the ${\delta}-$approximation of $f$, where
$f\in{R(A)}$.  Let $Au=f$ be solvable and let $y$ be its
minimal-norm solution.  The problem is: given $\{{\delta}, A,
f_{\delta}\}$, find the stable solution $u_{\delta}$ such that
error estimate (\ref{E4:unois}) holds.  Let
\begin{equation}\label{E36:b}
Bu=q:=A^{*}f, \quad where \quad B=A^{*}A\geq{0}.
\end{equation}
Let $q_{\delta}$ be given in place of $q$. Since $Au=f$ is
solvable, it is equivalent to $Bu=q$.  Since A is injective, B is
also injective. Assume without loss of generality $\|A\|\leq{1}$,
which implies that $\|A^{*}\|\leq{1}$.  Since
$\|f-f_{\delta}\|\leq{\delta}$ we obtain,
$\|q-q_{\delta}\|\leq\|A^{*}\|\delta$, hence we obtain
\begin{equation}\label{E38:norm}
\|q-q_{{\delta}}\|\leq{\delta}.
\end{equation}

Consider the iterative process:
\begin{equation}\label{E37:iter}
u_{n+1}=u_{n}-\mu{(Bu_{n}-q)}, \quad 0<\mu{<}\frac{1}{\|B\|},\quad
u(0)=u_{0}\perp{N(A)}.
\end{equation}
For example one may take $u_{0}=0$.  We obtain the following
result:
\begin{lemma}\label{l3:sol}
Assume that equation (\ref{E1:oper}) is solvable, and that $y$ is
its minimal-norm solution.  Then
\begin{equation}\label{E37:sol}
lim_{n\rightarrow{\infty}}{u_{n}}=y.
\end{equation}
\end{lemma}
\emph{Proof}

We note that from equation (\ref{E37:iter}),
\begin{equation}\label{E38:sl}
y=y-\mu{(By-q)}.
\end{equation}
Denote $u_{n}-y:=\gamma_{n}$.  Subtracting equation (\ref{E38:sl})
from equation (\ref{E37:iter}) and using induction,\\ we obtain

$\gamma_{n+1}=\gamma_{n}-\mu{B}\gamma_{n}=(I-\mu{B)}\gamma_{n}$
$=...=(I-\mu{B})^{n+1}\gamma_{0}$

with $\gamma_{0}=u_{0}-y,\quad\gamma_{0}\perp{N(A)}$.  Since,
$0<(1-\mu{\lambda})<1$, for all $\lambda{\in{(0,\|B\|)}}$, we
have,

$\|\gamma_{n}\|{^{2}}=\int_{0}^{\|B\|}{|1-\mu{\lambda}|^{2n}}d(E_{\lambda}\gamma_{0},
\gamma_{0})$
$=\int_{0}^{\epsilon}+\int_{\epsilon}^{\|{B}\|}:=I_{1}+I_{2}$.

If $\epsilon{\leq}\lambda{\leq}\|B\|$ then
$1-\mu{\lambda}\leq{1-\mu\epsilon}<1$.  Denote\\
$p:={1-\mu\epsilon}, 0<p<1.$\\ Then\\
$I_{2}\leq{p^{2n}\longrightarrow{0}},\quad
n\longrightarrow{\infty}$.\\ So\\
$I_{2}\longrightarrow{0},\quad{n\longrightarrow{\infty}}.$ \\
Since $1-\mu{\lambda}<1,$\\
$I_{1}\leq{\int_{0}^{\epsilon}}d(E_{\lambda}\gamma_{0},\gamma_{0})\longrightarrow{0},\quad{\epsilon\longrightarrow{0}}
$.\\  Since $\gamma_{0}\perp{N(B)=N(A)}$,
$\|{\gamma}_{n}\|\longrightarrow{0}$, as
$n\longrightarrow{\infty}$. Hence (\ref{E37:sol}) holds.  $\Box$

Now, we shall prove the main theorem,
\begin{theorem}\label{t2:iter}
Suppose that $A$ is a linear, bounded, injective operator on a
Hilbert space with $A^{-1}$ unbounded satisfying the equation
$Au=f$.  If $q_{\delta}$ is given such that
$\|q_{\delta}-q\|\leq{\delta},$ then one can use the iterative
process (\ref{E37:iter}) with $q_{\delta}$ in place of $q$ for
constructing a stable approximation of the solution $y$.
\end{theorem}
\emph{Proof:} By the iterative process, (\ref{E37:iter})
\begin{equation}\label{E39:rec}
u_{n+1,\delta}=u_{n,\delta}-\mu{(Bu_{n,\delta}-q_{\delta})},\quad
u(0)=u_{0}.
\end{equation}
From equation (\ref{E38:sl}) one has:\\ $y=y-\mu{(By-q)}$.\\
Denote $u_{n,\delta}-y:=\gamma_{n,\delta}.$\\  Then subtracting
equation (\ref{E38:sl}) from equation (\ref{E39:rec}),
\begin{equation}\label{E40:iter}
{\gamma}_{n+1,\delta}={\gamma}_{n,\delta}-\mu{B\gamma_{n,\delta}+\mu{(q_{\delta}-q)}},\quad\gamma_{0}=u_{0}-y.
\end{equation}
So that by induction,

$\gamma_{n,{\delta}}$
$={(I-\mu{B})}^{n}\gamma_{0}+\sum_{j=0}^{n-1}{(\mu{B})}^{j}\mu{(q_{\delta}-q)},\quad$
$\gamma_{0}=u_{0}-y$.

Thus,
\begin{equation}\label{E41:sum}
\gamma_{n,{\delta}}=\gamma_{n}+\sum_{j=0}^{n-1}{(\mu{B})}^{j}\mu{(q_{\delta}-q)},\quad{\gamma_{0}=u_{0}-y}.
\end{equation}

Since, $\|\mu{B}\|\leq{1}$ and by using (\ref{E38:norm}),  we
obtain,

$\|\gamma_{n,\delta}\|\leq\|\gamma_{n}\|+n\mu{\delta},\quad{n}\geq{1}$

It is already proved in lemma (\ref{l3:sol}), that
$\|\gamma_{n}\|\longrightarrow{0}\quad{as}\quad{n}\longrightarrow{\infty}$.
Hence $\|\gamma_{n,\delta}\|\longrightarrow{0}$ as
$\delta{\longrightarrow{0}}.$  Thus theorem \ref{t2:iter} is
proved.  $\Box$

\emph{Remark}:  In this method the regularization parameter is the
stopping rule, $n({\delta})$, the number of iterations and can be
found by solving the minimization problem

$\|\gamma_{(n)}\|+n\mu{\delta}=min.\longrightarrow{0}$ as
$\delta\longrightarrow{0},\quad{n\geq{1}}$ and
$n({\delta})\longrightarrow{\infty}$ as
$\delta\longrightarrow{0}$.
\newpage
\section{Dynamical systems method}
In this section we study dynamical systems method for solving
linear and non-linear ill-posed problems in a real Hilbert space
H. The DSM for solving operator equations consists of a
construction of a Cauchy problem, which has a unique global
solution for an arbitrary initial data, this solution tends to a
limit as time tends to infinity, and this limit is the stable
solution of the given operator equation.  This method can be used
for solving well-posed problems also.  Our discussion is based on
the paper by A.G.Ramm \cite{aR04}.

Consider an operator equation
\begin{equation}\label{E42:dsm}
F(u):=Bu-f=0,\quad{f}\in{H}
\end{equation}
where $B$ is a linear or non-linear operator in a real Hilbert
space $H$.  We make the following assumptions.

\emph{Assumption 1} Assume that $F$ has two Fr\'echet
$u_{\alpha{,\delta}}$ derivatives: $F\in{C_{loc}^{2}}$, i.e.,
\begin{equation}\label{E43:fre}
\sup_{u\in{B(u_{0},R)}}\|F^{(j)}(u)\|\leq{M_{j}(R)}, \quad j=0,1,2
\end{equation}
where $B(u_{0},R):=\{u:\|u-u_{0}\|\leq{R}\},\quad{u_{0}}$ is
arbitrary fixed element in $H$ and $R>0$ is arbitrary and
$F^{(j)}(u)$ is the j-th Fr\'echet derivative of $F(u).$

\emph{Assumption 2} Assume that there exists a  solution
$y\in{B(u_{0},R)}$ (not necessarily unique globally) to equation
(\ref{E42:dsm}):
\begin{equation}\label{E44:exact}
F(y)=0
\end{equation}

Problem (\ref{E42:dsm}) is called \emph{well-posed} if $F'(u)$ is
a bounded invertible linear operator,  i.e., if $[F'(u)]^{-1}$
exists and if the estimate
\begin{equation}\label{E45:well}
\sup_{u\in{B(u_{0},R)}}\|{[F^{'}(u)]}^{-1}\|\leq{m(R)},
\end{equation}
Otherwise, it is called \emph{ill-posed.}

Let $\dot{u}$ denote time-derivative.  Consider the Cauchy problem (dynamical system):
\begin{equation}\label{E46:cauchy}
\dot{u}=\mathbf{\Phi}(t,u);\quad{u(0)=u_{0}}
\end{equation}
where $\mathbf{\Phi}$ is a non-linear operator, which is locally
Lipschitz with respect to $u\in{H}$ and continuous with respect to
$t\geq{0}$, so that the Cauchy problem (\ref{E46:cauchy}) has a
unique local solution.  The operator $\mathbf{\Phi}$ is chosen
such that the following properties hold:
\begin{enumerate}
\item There exists unique global solution $u(t)$ to
the Cauchy problem (\ref{E46:cauchy}).  (Here global solution
means the solution defined for all $t>0$. )
\item There exists $u(\infty{)}:=\lim_{t\rightarrow{\infty}}{u(t)}.$
\item and finally this limit solves equation (\ref{E42:dsm}):
$F(u(\infty{)})=0.$
\end{enumerate}
Problem (\ref{E42:dsm}) with noisy data
$f_{\delta},\quad{\|f_{\delta}-f\|}\leq\delta$, given in place of
$f$, generates the problem:
\begin{equation}\label{E48:prob}
\dot{u}_{\delta}= \mathbf{\Phi}_{\delta}{(t,u_{\delta})},\quad{u_{\delta}(0)}=u_{0},
\end{equation}
The solution $u_{\delta}$ to problem (\ref{E48:prob}), calculated
at $t=t_{\delta}$, where $t_{\delta}$ is suitably chosen,
satisfies the error estimate
\begin{equation}\label{E49:lim}
\lim_{{\delta}\to{0}}\|u_{\delta}(t_{\delta})-y\|=0.
\end{equation}
The choice of $t_{\delta}$ with this property is called the
\emph{stopping rule} and is the regularization parameter in DSM
method.  One has usually
$\lim_{{\delta}\to{0}}t_{\delta}={\infty}$.

Dynamical systems method can be used to solve ill-posed and also
well-posed problems.  In this report we are interested in
discussing solving linear ill-posed problems by DSM.  One can also
find in A.G.Ramm's paper \cite{aR04}, a discussion of DSM for
solving well-posed problems, nonlinear ill-posed problems with
monotone and non-monotone operators and the recent development of
the theory of DSM.

\chapter{Dynamical systems method for linear problems}

In this section, for linear solvable ill-posed problem $Au=f$,
with bounded linear operator $\|A\|<1$, DSM is justified and a
stable approximation of the minimal norm solution to ill-posed
problem  with noisy data $f_{\delta}$,
$\|f_{\delta}-f\|\leq{\delta}$ is constructed. This section is
based on paper \cite{aR04}.

Assume that (2.24) and (2.25) holds and (2.26) fails so the
problem is ill-posed.  Consider the equation
\begin{equation}\label{E51:op}
Au=f
\end{equation}
where $f\in{R(A)}$ is arbitrary.  Let us assume the following

\emph{Assumptions}(*).
\begin{enumerate}
\item Let $A$ be a linear, bounded operator in a Hilbert space $H$, defined on all of $H$, the  range $R(A)$ is not closed, so that $A^{-1}$ is unbounded.  So problem (\ref{E51:op}) is an ill-posed
problem. Let $f_{\delta}$ be given in place of $f$,
$\|f-f_{\delta}\|\leq{\delta}$.
\item Equation (\ref{E51:op}) is solvable (possibly non-uniquely). Let $y$ be the  minimal-norm solution to equation (\ref{E51:op}), ${y}\perp{N(A)}$, where $N(A):=\{v:Av=0\}$ is the null-space of $A$.
\end{enumerate}
Let $B=A^{*}A\geq{0}$ and $q:=A^{*}f, {A}^{*}$ is the adjoint of
$A$. Then we obtain the normal equation,
\begin{equation}\label{E52:eqi}
Bu=q.
\end{equation}
We know that if equation (\ref{E51:op}) is solvable then it is
equivalent to equation (\ref{E52:eqi}) with $q_{\delta}$ given in
place of $q$. Without loss of generality assume $\|A\|\leq{1}$, so
$\|A^{*}\|\leq{1}$ and $\|B\|\leq{1}$.  Then
$\|q-q_{\delta}\|=\|A^{*}(f-f_{\delta})\|\leq\|A^{*}\|\delta\leq{\delta}$
and $y\perp{N(B)}$.  Let $\epsilon{(t)}>0$, be a continuous,
monotonically decaying function decaying to zero function on
$\mathbf{R}_{+}$ such that $\int_{0}^{\infty}\epsilon{ds}=\infty$.
Let
\begin{equation}\label{E53:newF}
F(u):=Bu-q=0
\end{equation}
then  $F^{'}(u)=B$.  Consider the Cauchy problem
\begin{equation}\label{E54:cau}
\dot{u}=\mathbf{\Phi}(t,u),\quad u(0)=u_{0}.
\end{equation}
$\mathbf{\Phi(t,u)}$
$:=-[F^{'}(u)+\epsilon{(t)}]^{-1}[F(u)+\epsilon{(t)}u]$
$=-[B+\epsilon{(t)}]^{-1}[Bu-q+\epsilon{(t)}u]$.

Thus from the Cauchy problem (\ref{E54:cau}), the DSM for solving
equation (\ref{E53:newF}) is solving the Cauchy problem
\begin{equation}\label{E54:cuchy}
\dot{u}=\mathbf{\Phi}(t,u)=-u+[B+\epsilon{(t)}]^{-1}q,\quad{u(0)=u_{0}}.
\end{equation}
with
\begin{equation}\label{E55:dela}
\mathbf{\Phi}_{\delta}=-u_{\delta}+[B+\epsilon{(t)}]^{-1}q_{\delta}.
\end{equation}

We now prove the main theorem of this section: \emph{Given
noisy data $f_{\delta}$, every linear ill-posed problem
(\ref{E51:op}) under the assumptions (*) can be
stably solved by the DSM.}

\begin{theorem}\label{t1:lin}
Assume (*), and let $B:=A^{*}A$, $q:=A^{*}f$.  Assume
$\epsilon{(t)>0}$ to be a continuous, monotonically decaying to
zero function on $[0,{\infty})$ such that
$\int_{0}^{\infty}\epsilon{ds}=\infty$.  Then we have the
following results.
\begin{enumerate}
\item For any $u_{0}\in{H},$ the Cauchy problem (\ref{E54:cuchy}) has a unique global solution $u(t)$, (the initial approximation $u_{0}$ need not be close to the solution $u(t)$ in any sense).
\item There exists $\lim_{t\rightarrow{\infty}}u(t)=u({\infty})=y$, and $y$ is the unique minimal-norm solution to equation (\ref{E51:op}).  ${Ay=f},y\perp{N}$, and $\|y\|\leq\|z\|$, for all $z\in{N:=\{{z:F(z)=0\}}}$.
\item If $f_{\delta}$ is given in place of $f$, $\|f-f_{\delta}\|\leq{\delta}$,
then there exists a unique global solution $u_{\delta}(t)$ to the
Cauchy problem
\begin{equation}\label{56:linedel}
\dot{u}_{\delta}=\mathbf{\Phi}_{\delta}(t,u_{\delta})=-u_{\delta}+[B+\epsilon{(t)}]^{-1}q_{\delta},
\quad{u_{\delta}}(0)=u_{0}
\end{equation}
with $q_{\delta}:=A^{*}f_{\delta}$.
\item There exists $t_{\delta}$, such that it satisfies the error estimate
\begin{equation}\label{57:udel}
\lim_{\delta\rightarrow{0}}\|u_{\delta}(t_{\delta})-y\|=0,\quad{\lim_{\delta\rightarrow{0}}}t_{\delta}={\infty}.
\end{equation}
The choice of $t_{\delta}$ with this property is the stopping
rule.  This $t_{\delta}$ can be for example chosen  by a
discrepancy principle or as a root of the equation
\begin{equation}\label{E58:tdel}
2\surd{\epsilon{(t)}}={\delta}^{b},\quad{b}\in(0,1).
\end{equation}
\end{enumerate}
\end{theorem}
\emph{Proof}: Since the Cauchy problem (\ref{E54:cau}) is linear,
its solution can be written by an explicit analytic formula
\begin{equation}\label{E59:intsol}
u(t)=u_{o}e^{-t}+\int_{0}^{t}e^{-(t-s)}[B+\epsilon{(s)}]^{-1}Byds.
\end{equation}

Taking limit as $t\rightarrow{\infty}$ to (\ref{E59:intsol}) and
applying L'Hospital's rule, to the second term  in the right hand
side of equation (\ref{E59:intsol}), we obtain,

$\lim_{t\rightarrow{\infty}}\frac{\int_{0}^{t}e^{s}[B+\epsilon{(s)}]^{-1}Byds}{e^{t}}$
$=\lim_{t\rightarrow{\infty}}[B+\epsilon{(t)}]^{-1}By,$

provided only that $\epsilon{(t)}>0$ and
$\lim_{t\rightarrow{\infty}}\epsilon{(t)}=0$.\\ Since
$y\perp{N}=N(B)=N(A)$,
\begin{equation}\label{E60:lop}
\lim_{{\epsilon}\rightarrow{0}}{[B+\epsilon{]}}^{-1}By
=\lim_{\epsilon{\rightarrow{0}}}\int_{0}^{\|B\|}\frac{\lambda}{{\lambda}+\epsilon}dE_{\lambda}y
=\int_{0}^{\|B\|}dE_{\lambda}y
=y,
\end{equation}
by the spectral theorem and by the Lebesgue dominated convergence
theorem, where $E_{\lambda}$ is the resolution of the identity
corresponding to the self-adjoint operator $B$, $\lambda$  is
taken over the spectrum of $B$, and
${\lim}_{\epsilon{\rightarrow}{0}}\frac{\lambda}{\lambda{+}\epsilon}=1,$
for ${\lambda}>0$ and $=0,$ for ${\lambda}=0$.  Thus from
equations (\ref{E59:intsol}) and (\ref{E60:lop}) there exists
$u({\infty})=\lim_{t\to\infty}u(t)=y$ with $Ay=f$.

Denote $\eta{(t)}:=\|u(t)-y\|$, then
$\lim_{t\rightarrow{\infty}}\eta{(t)}=0$.  In general, the rate of
convergence of ${\eta}$ to zero can be arbitrarily slow for a
suitably chosen $f$.  Under an additional a priori assumption of
$f$ (for example, the source type assumptions), this rate can be
estimated.

\emph{Proof of results 3 and 4:  Derivation of the stopping rule.}

Consider the Cauchy problem with noisy data.  Suppose $f_{\delta}$
is given, with $\|f_{\delta}-f\|\leq{\delta}$, then
$\|q_{{\delta}}-q\|\leq\delta$.  We require the following lemma
for the proof:
\begin{lemma}\label{l6:lem}
\begin{equation}\label{E65:lemm}
\|[B+{\epsilon}]^{-1}A^{*}\|\leq{\frac{1}{2\surd{\epsilon}}}
\end{equation}
\end{lemma}
\emph{Proof of lemma}:  We have $B:=A^{*}A\geq{0}$.  Denote
$Q:=AA^{*}\geq{0}$, then $\|Q\|\leq{1}$.   We have

$[B+\epsilon{I}]^{-1}A^{*}=A^{*}[Q+\epsilon{I}]^{-1}.$

So,

$\|[B+\epsilon{I}]^{-1}A^{*}\|=\|A^{*}[Q+\epsilon{I}]^{-1}\|=\|UQ^{1/2}[Q+\epsilon{I}]^{-1}\|$,

(U being the isometry operator $(\|Uf\|=\|f\|)$, and
$A=U|A|=U(A^{*}A)^{1/2}$,polar representation of the linear
operator).

$=\|Q^{1/2}[Q+\epsilon{I}]^{-1}\|$, (since $\|U\|=1$.)
$=\|\int_{0}^{1}{\frac{{\lambda}^{1/2}}{\lambda{+\epsilon}}dE_{\lambda}\|}$
$=\sup_{0<{\lambda}\leq{1}}\frac{{\lambda}^{1/2}}{\lambda{+\epsilon}}$
$=\frac{1}{2\surd{\epsilon}}.$

(by the spectral theorem:
$\phi{(A)}=\int_{-\infty}^{\infty}\phi{({\lambda})}dE_{\lambda},$
$\|\phi{(A)}\|=\sup_{\lambda}|\phi{(\lambda{)}}|$, ${\lambda}$ is
taken over the spectrum of $A$). $\Box$

By triangle inequality,

$\|u_{\delta}(t)-y\|$ $=\|u_{\delta}(t)-u(t)+u(t)-y\|$
$\leq\|u_{\delta}(t)-u(t)\|+\|u(t)-y\|$
$=\|u_{\delta}(t)-u(t)\|+\eta{(t)}.$

$\|u_{\delta}(t)-u(t)\|$
$=\|\int_{0}^{t}e^{-(t-s)}[B+\epsilon{(s)}]^{-1}(q_{\delta}-q)ds\|$
$\leq\int_{0}^{t}e^{-(t-s)}\frac{\delta}{2\surd{\epsilon}}$, by
(\ref{E65:lemm}) $\leq\frac{\delta}{2\surd{\epsilon}}$, since
$\int_{0}^{t}e^{-(t-s)}ds=1-e^{-t}\leq{1}.$

So,
\begin{equation}\label{E66:dela}
\|u_{\delta}(t)-u(t)\|\leq{\frac{\delta}{2\surd{\epsilon}}}.
\end{equation}
Thus,
$\|u_{\delta}(t)-y\|\leq{\frac{\delta}{{2\surd{\epsilon}}}}+\eta{(t)}$.

We have already proved that
$\lim_{{t\rightarrow{\infty}}}\eta{(t)}=0$.  Choose,
$t=t_{\delta}$, satisfying equation (\ref{E58:tdel}) then this
particular choice of $t_{\delta}$ satisfies the error estimate
(\ref{57:udel}) and
$\|u_{\delta}(t)-u(t)\|\longrightarrow{0}\quad{as}\quad\delta\longrightarrow{0}$.
If the decay rate of $\eta{(t)}$ is known, a more efficient
stopping rule can be obtained by choosing, $t=t_{\delta}$ such
that it satisfies the minimization problem

$\frac{\delta}{2\surd{\epsilon}}+\eta{(t_{\delta})}=$ min.
$\longrightarrow{0}$ as ${\delta}\longrightarrow{0}$. $\Box$

$\mathbf{Remarks}$

Remark 1: \emph{Discrepancy principle for the DSM}.

Choose the stopping time $t_{\delta}$ as the  unique solution to
the equation: $\|Au_{\delta}(t)-f_{\delta}\|=C{\delta}$ where
C=constant$>1$, where it is assumed that
$\|f_{\delta}\|>{\delta}$.  In-addition, we assume that
$f_{\delta}\perp{N(A^{*})}$, so that C=1, then the equation is:
\begin{equation}\label{E66:disp}
\|A[B+\epsilon{(t)}]^{-1}A^{*}f_{\delta}-f_{\delta}\|={\delta}.
\end{equation}
Then this $t_{\delta}$ satisfies the error estimate
(\ref{57:udel}). One can find detailed discussion of this in
A.G.Ramm's paper \cite{aR04a}.

Remark 2: \emph{Choosing scaling parameter $\epsilon{(t)}$}.

We can choose the scaling parameter as large as we wish. In
particular we can choose
\begin{equation}\label{E67:epsi}
\epsilon{(t)}=\frac{c_{1}}{(c_{0}+t)^{b}}
\end{equation}
where, $0<b<1, c_{1}, c_{2}$ are positive constants.

\backmatter \appendix


\begin{thebibliography}{9}
\bibitem{aR02}
Ramm,A.G., \emph{Regularization of ill-posed problems with
unbounded operators,} J.Math.Anal.Appl. 271, (2002), 547-550.
\bibitem{aR04}
Ramm,A.G., \emph{Dynamical systems method for solving operator
equations,} Communic. in Nonlinear Sci. and Numer. Simulation, 9,
N2, (2004), 547-550.
\bibitem{aR03}
Ramm,A.G., \emph{On a new notion of regularizer,} J.Phys.A. 36,
(2003), 2191-2195.
\bibitem{aR03a}
Ramm,A.G.,
\emph{Global convergence for ill-posed equations with monotone operators: the dynamical systems method,}
 J.Phys A, 36, (2003), L249-L254.
\bibitem{aR04a}
Ramm,A.G., \emph{Lecture Notes,Math 994, Ill-Posed And Inverse
Problems,} Department of Mathematics, Kansas State University,
(Fall 2003, Spring 2004).
\bibitem{aR90}
Ramm,A.G., \emph{Random Fields Estimation Theory,} Longman/Wiley,
NY,(1990), 106-134.
\bibitem{aR03a1}
Ramm,A.G., \emph{On the discrepancy principle,} Nonlinear Funct.
Anal. Appl., 8, N2, (2003), 307-312.
\bibitem{aM84}
Morozov,V.A., \emph{Methods For Solving Incorrectly Posed
Problems,} Springer Verlag, NY, (1984).
\bibitem{aKrT94}
Kaplan.A. and Tichatschke.R., \emph{Stable Methods For Ill-Posed
Variational Problems,} Akademie Verlag, Berlin, (1994).
\bibitem{mLgRpS86}
Lavrent'ev.M.M., Romanov.V.G., and Shishat.skii.S.P.,
\emph{Ill-Posed problems Of Mathematical Physics and Analysis,}
Translations Of Mathematical Monographs, 64, AMS., Providence,
Rhode Island (1986).
\bibitem{rGkH83}
Gorenflo.R. and Hoffmann.K.H., \emph{Applied Nonlinear Functional
Analysis, Variational Methods And Ill-Posed Problems,} VerlagPeter
lang, Frankfurt, (1983).
\bibitem{gHkH83}
Hammerlin.G. and Hoffmann.K.H., \emph{Improperly Posed Problems
And Their Numerical Treatment,} International Series Of Numerical
Mathematics, 63, Birkhauser,(1983).
\bibitem{aR82}
Ramm,A.G., \emph{Iterative Methods For Calculating Static Fields
And Wave Scattering By Small Bodies,} Springer Verlag, NY.,
(1982).
\bibitem{aRaK}
Ramm.A.G. and Katsevich.A.I., \emph{The Radon Transform And Local
Tomography,} CRC., NY., (1996).
\bibitem{hEcG87}
Engl.H.W. and Groetsch.C.W., \emph{Inverse And Ill-Posed
Problems,} Notes And Reports In Mathematics In Science And
Engineering, 4, Academic Press, Boston (1987).
\bibitem{aR04a1}
Ramm.A.G., \emph{Discrepancy principle for the dynamical systems
method}, Communic. in Nonlinear Sci. and Numer.Simulation, (to
appear).
\bibitem{aR68}
Ramm.A.G, \emph{On numerical differentiation},
Mathematics,Izvestija VUZov, 11, (1968), 131-135.\end{thebibliography}
\end{document}